\documentclass{amsart}

\usepackage{amsmath,amssymb,graphicx,color,float,tikz,latexsym}

\newtheorem{Theorem}{Theorem}[section]
\newtheorem{Lemma}[Theorem]{Lemma}
\newtheorem{Corollary}[Theorem]{Corollary}
\newtheorem{Proposition}[Theorem]{Proposition}
\newtheorem{Remark}[Theorem]{Remark}

\newtheorem{Definition}[Theorem]{Definition}
\newtheorem{Problem}[Theorem]{Problem}

\newenvironment{pf}{\medskip\par\noindent{\bf Proof\/}.}{\hfill
$\Box$}

\begin{document}
\begin{center}
\renewcommand{\thefootnote}{\fnsymbol{footnote}}
{\large\bf Projective dimension and  regularity  of the path ideal of the line graph}
\footnote[2]{Supported by the National Natural Science Foundation
of China (11271275) and  by  Foundation of Jiangsu Overseas Research \& Training Program for
University Prominent Young \& Middle-aged Teachers and Presidents and
 by Foundation of the Priority Academic Program Development of Jiangsu Higher Education Institutions. }
\\
\vspace{1.0cm}

\renewcommand{\thefootnote}{}
 Guangjun Zhu \footnote{ School of Mathematical Sciences, Soochow University,
 Suzhou 215006, P. R. China,
e-mail:
zhuguangjun@suda.edu.cn} \\

\end{center}

\vspace{0.5cm}

{\bf Abstract:} By generalizing the notion of the path ideal of a graph,
we study  some algebraic properties of some path ideals associated to a line graph.
We  show that the quotient ring of these ideals are always sequentially Cohen-Macaulay and
 also provide some exact formulas for the projective dimension and the regularity
of these ideals. As some consequences, we give  some exact formulas for the depth
of these ideals.

\vspace{0.5cm}


\vspace{0.5cm}


\vspace{0.5cm}

\section{Introduction }

\hspace{3mm} The path ideal of a graph was first introduced by Conca and De Negri \cite{CN}.
Fix an integer $m\geq 2$, and suppose that $\Gamma$ is a directed graph with vertex set $V=\{x_{1},\dots,x_{n}\}$, i.e., each edge has been
assigned a direction. A sequence of $m$ vertices $x_{i_{1}},\dots,x_{i_{m}}$
is said to be a path of length $m$  if there are $m-1$  distinct edges $e_{1},\dots,e_{m-1}$ such that $e_{j}=(x_{i_{j}},x_{i_{j+1}})$
is a directed
edge from $x_{i_{j}}$
to $x_{i_{j+1}}$. By identifying the vertices with the variables in the polynomial
ring $R=k[x_{1},\dots, x_{n}]$	 over a field $k$, the path ideal of $\Gamma$  of length $m$ is the monomial ideal
$$J_{m}(\Gamma)=
(\{x_{i_{1}}\cdots x_{i_{m}}\ |\ x_{i_{1}},\dots ,x_{i_{m}}
\mbox{is a path of length}\ m\  \mbox{in}\ \Gamma
\})$$
Note that when $m=2$, then $J_{2}(\Gamma)$ is simply the edge ideal of $\Gamma$, which is defined by Villarreal in \cite{V1}.
Other   higher dimensional analogues can be found in \cite{F2,HT2}, among others. The underlying theme in all correspondences is to relate the algebraic properties
to the combinatorial properties, and vice versa. We mainly
study the algebraic properties of the path ideal.

Path ideals appeared in \cite{CN} as an example of a family of monomial ideals that
are generated by $M$-sequences. Among other things, it is shown that when $\Gamma$ is a
directed tree, the Rees algebra $\mathcal{R}(J_{m}(\Gamma))$ is normal and Cohen-Macaulay.
The path ideals of complete bipartite graphs are shown to be normal in \cite{RV}, while the path
ideals of cycles are shown to have linear type in \cite{BS}. In \cite {HT3}, He and  Tuyl
study $J_{m}(\Gamma)$ in the special case that  $\Gamma$  is the line graph $L_{n}$.
The line graph $L_{n}$ is a graph with vertex set $V=\{x_{1},\dots,x_{n}\}$ and directed edges
$e_{j}=(x_{j},x_{j+1})$ for $j=1,\dots,n-1$. Thus, the graph $L_{n}$ has the form

\begin{center}
\begin{figure}[!h]
\vspace{-0.5cm}
\includegraphics[height=2cm,width=10cm]{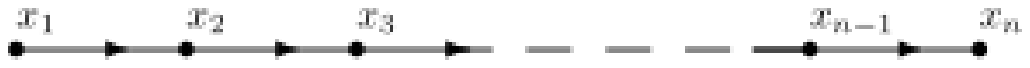}\\
\vspace{-0.5cm}
\end{figure}
\end{center}

\vspace{-0.5cm}They prove that $R/J_{m}(L_{n})$ is sequentially Cohen-Macaulay and also provide an exact formula for the projective dimension
of $J_{m}(L_{n})$ in terms of $m$ and $n$. They showed that:

 \begin{Theorem}(Theorem 4.1)
 Let $p,m,n,d$
be  integers such that  $n=p(m+1)+d$, where $p\geq 0$, $0\leq d\leq m$ and $2\leq m\leq n$. Then
the projective dimension of $J_{m}(L_{n})$ is given by
$$\mbox{pd}\,(J_{m}(L_{n}))=\left\{\begin{array}{ll}
2p-1& d\neq m;\\
2p& d=m.
\end{array}\right.$$
\end{Theorem}

\vspace{3mm}In \cite{AF}, using purely combinatorial arguments, Alilooee and Faridi also gave the above formula for
 projective dimension of $J_{m}(L_{n})$. Furthermore,  they  gave an explicit formula for
Castelnuovo-Mumford regularity  of $J_{m}(L_{n})$ in terms of $m$ and $n$. They showed that:

 \begin{Theorem}(Corollary 4.14)
 Let $p,m,n,d$
be  integers such that  $n=p(m+1)+d$, where $p\geq 0$, $0\leq d\leq m$ and $2\leq m\leq n$. Then
the regularity of $J_{m}(L_{n})$ is given by
$$\mbox{reg}\,(J_{m}(L_{n}))=\left\{\begin{array}{ll}
p(m-1)+1& d\neq m;\\
p(m-1)+m& d=m.
\end{array}\right.$$
\end{Theorem}

\vspace{3mm}We generalize the notion of the path ideal as the following:
Let $\Gamma$ be  a directed graph
with vertex set $V=\{x_{1},\dots,x_{n}\}$,  the path ideal of $\Gamma$  of length $m$ is the monomial ideal
$$I_{m,k}(\Gamma)=(u_{1},\dots,u_{k}),\  \mbox{where}\  u_{1},\dots,u_{k}
\ \mbox{ are some paths of length}\  m\ \mbox{in}\ \Gamma.$$
When $u_{1},\dots,u_{k}$ are all paths of length $m$ in $\Gamma$, $I_{m,k}(\Gamma)=J_{m}(\Gamma)$.

To the best of our knowledge,
little  is known about these ideals.  It is, therefore, of interest to determine
algebraic properties of the ideals $I_{m,k}(\Gamma)$. In this article we shall focus on the case
that $\Gamma$  is the line graph $L_{n}$ and $I_{m,k}(L_{n})=(u_{1},\dots,u_{k})$,  where for any $1\leq i\leq k$, $u_{i}=\prod\limits_{j=1}^{m}x_{(i-1)(m-l)+j}$
 is a path of length $m$ in $L_{n}$  and $1\leq l\leq m$ is an  integer.  we shall abuse notation
and write  $I_{m,k}(L_{n})$ for $I_{m,l,k}$.
In Section 2, we study algebraic properties of the ideal  $I_{m,l,k}$ and
show that $R/I_{m,l,k}$ is sequentially Cohen-Macaulay. In Section 3, using the notion of a Betti-splitting, as defined
in \cite {FHT}, we  derive some exact formulas for the projective dimension  and  regularity  of the ideal $I_{m,l,k}$ (see Theorems \ref{Thm3}, \ref{Thm4} and \ref{Thm5}).
  As some consequences, we give  some exact formulas for the depth
of these ideals.

\vspace{0.5cm}

\section{Preliminaries}

\hspace{3mm}In this section, we will show that the ideal $I_{m,l,k}$
can be viewed as the facet ideal of the simplicial complex $\Delta_{m,l,k}$ or the edge ideal of the clutter $\mathcal{C}_{m,l,k}$. By proving $\mathcal{C}_{m,l,k}$ has  the free vertex property,  we can obtain that the quotient ring $R/I_{m,l,k}$ is sequentially Cohen-Macaulay. We recall the relevant definitions.

\vspace{3mm}\begin{Definition}
A simplicial complex $\Delta$ on the vertex set $V$ is a collection of subsets of  $V$  with the property that if $F\in \Delta$
then all subsets of $F$ are also in $\Delta$. An element of $\Delta$ is called a face, the dimension of a face $F$ is $|F|-1$,
and the dimension of $\Delta$ is the largest dimension of  faces  of $\Delta$. The   maximal faces of $\Delta$ under inclusion
are called facets, and the set of facets of $\Delta$ is denoted by Facets\,$(\Delta)$. Simplicial complex $\Delta$ is called
pure if all of its facets have the same dimension, otherwise
$\Delta$ is non-pure. If Facets\,$(\Delta)=\{F_{1},\dots,F_{q}\}$
 we write $\Delta=\langle F_{1},\dots,F_{q}\rangle$.
\end{Definition}

\vspace{3mm}\begin{Definition}
A clutter $\mathcal{C}$ on vertex set $V$ is a family of subsets of $V$, called edges,
none of which is included in another. The set of vertices and edges of $\mathcal{C}$  are denoted by
$V_{\mathcal{C}}$  and $E_{\mathcal{C}}$  respectively.
\end{Definition}

\vspace{3mm}Given a clutter $\mathcal{C}=(V_{\mathcal{C}},E_{\mathcal{C}})$, we can associate to $\mathcal{C}$ the simplicial complex
$\Delta=\{F\subseteq V_{\mathcal{C}}\mid F\subseteq E_{i}\ \text{for some}\ \ E_{i}\in E_{\mathcal{C}}\}$.
Conversely, given a simplicial complex $\Delta$ with vertex set $V$ and set of facets Facets\,$(\Delta)$,
we can associate to $\Delta$  the clutter $\mathcal{C}=(V,\text{Facets}\,(\Delta))$.

\vspace{3mm}Let  $I$ be a squarefree monomial ideal of $R$
with minimal generators $x^{v_{1}},\dots,x^{v_{q}}$. We
use $x^{a}$ as an abbreviation for $x_{1}^{a_{1}}\cdots x_{n}^{a_{n}}$, where $a=(a_{1},\dots, a_{n})\in \mathbf{N}^{n}$. Note that the
entries of each $v_{i}$ are in $\{0,1\}$. We associate to the ideal $I$ a clutter  $\mathcal{C}$ by taking the set
of indeterminates  $V_{\mathcal{C}}=\{x_{1},\dots, x_{n}\}$ as the vertex set and $E_{\mathcal{C}}=\{S_{1},\dots, S_{q}\}$ as the edge
set, where $S_{i}=\text{supp}\,(x^{v_{i}})$ is the support of $x^{v_{i}}$, i.e., $S_{i}$ is the set of variables that occur in
$x^{v_{i}}$. For this reason $I$ is called the edge ideal of $\mathcal{C}$ and is denoted $I=I(\mathcal{C})$. Edge ideal
of a clutter is also called facet ideal   because $\{S_{1},\dots, S_{q}\}$ is exactly the  set of facets of the
simplicial complex $\Delta=\langle S_{1},\dots, S_{q}\rangle$ generated by $S_{1},\dots, S_{q}$.

\vspace{3mm}Let $\Delta$ be a simplicial complex and   $\sigma \in \Delta$,
  the deletion of $\sigma$ from $\Delta$ is the
simplicial complex defined by
$\Delta\setminus  \sigma=\{\tau\in \Delta|\  \sigma\not\subseteq \tau\},$
when  $\sigma=\{x\}$, we shall abuse notation and write
$\Delta\setminus x$ for $\Delta\setminus\{x\}$.
If $\Delta=\langle F_{1},\dots,F_{q}\rangle$, the simplicial complex obtained by
removing the facet $F_{i}$ from  $\Delta$ is the simplicial
complex
$\Delta\setminus \langle F_{i}\rangle=\langle F_{1},\dots,\hat{F_{i}},\dots,F_{q}\rangle.$

\vspace{3mm}The following definition of shellable is due to Bj\"orner and Wachs \cite{BW} and is
usually referred to as nonpure shellable, here we  drop the adjective
``nonpure".

\begin{Definition}  A simplicial complex
$\Delta$ is shellable if the facets  of $\Delta$ can be ordered
$F_{1},\dots,F_{s}$ such that for all $1\leq i<j\leq s$, there
exists some $x\in F_{j}\setminus F_{i}$ and some $l\in
\{1,\dots,j-1\}$ with $F_{j}\setminus F_{l}=\{x\}$. We call
$F_{1},\dots,F_{s}$  a shelling of $\Delta$ when the facets have
been ordered with respect to the shellable definition.
\end{Definition}

\vspace{3mm}
If the simplicial complex $\Delta$ is   pure
and  satisfies the above definition of shellable, we will say $\Delta$  is pure shellable.

\vspace{3mm}
\begin{Definition}\label{def9}
Let $M$ be a graded $R$-module.  $M$ is called sequentially
Cohen-Macaulay if there exists a filtration of
graded $R$-submodules of $M$
$$0=M_{0}\subset M_{1}\subset \cdots \subset M_{r}=M$$
such that each quotient $M_{i}/M_{i-1}$ is Cohen-Macaulay and the
Krull dimensions of the quotients are increasing, i.e.,
$\text{dim}\,(M_{1}/M_{0})<\text{dim}\,(M_{2}/M_{1})<\cdots<
\text{dim}\,(M_{r}/M_{r-1})$.

\vspace{0.3cm}
A simplicial complex is said to be sequentially Cohen-Macaulay if its
Stanley-Reisner ideal has a sequentially Cohen-Macaulay quotient.
\end{Definition}

\vspace{0.3cm}It is well known that   shellable implies
sequentially Cohen-Macaulay.

\vspace{0.3cm}
Let $\mathcal{C}$ be a clutter with vertex set $V$. A vertex cover of $\mathcal{C}$ is a subset $A$ of $V$ that intersects
every edge of $\mathcal{C}$. If $A$ is a minimal element (under inclusion) of the set of vertex covers of
$\mathcal{C}$, it is called a minimal vertex cover.
To a squarefree monomial ideal $I=I(\mathcal{C})$, it also corresponds to a simplicial
complex via the Stanley-Reisner correspondence \cite{V}. We let $\Delta_{\mathcal{C}}$ represent this simplicial
complex. Note that $F$ is a facet of $\Delta_{\mathcal{C}}$ if and only if $X\setminus F$ is a minimal vertex cover of
$\mathcal{C}$. As for clutters, we may say that the clutter $\mathcal{C}$ is shellable if $\Delta_{\mathcal{C}}$ is shellable.

\vspace{3mm}\begin{Definition}\label{def2} Let $I'\subsetneq I$ be two ideals of $R$,  $I'$  is called a minor of $I$ if there is a subset $V'=\{x_{i_{1}},\dots,x_{i_{r}},x_{j_{1}},\dots,x_{j_{s}}\}$
of the set of variables $V=\{x_{1},\dots, x_{n}\}$ such that $I'$ is a proper ideal of $R'=k[V\setminus V']$
that can be obtained from the generator set of $I$ by setting $x_{i_{k}}=0$ and $x_{j_{l}}=1$  for all $k,l$.
A minor of $\mathcal{C}$ is a clutter $\mathcal{C'}$ on the vertex
set $V_{\mathcal{C'}}= V\setminus V'$ that corresponds to a minor $(0)\subsetneq I'\subsetneq R'$. The edges of $\mathcal{C'}$
are obtained from $I'$ by considering the unique set of squarefree monomials of $R'$ that
minimally generate $I'$. For use below we say $x_{i}$ is a free variable (resp. free vertex) of $I$
(resp. $\mathcal{C}$) if $x_{i}$ only appears in one of the monomials $x^{v_{1}},\dots, x^{v_{q}}$ (resp. in one of the edges
of $\mathcal{C}$). If all the minors of $\mathcal{C}$ have free vertices, we say that $\mathcal{C}$ has the free vertex property.
Note that if $\mathcal{C}$ has the free vertex property, then so do all of its minors.
\end{Definition}

\vspace{3mm}Tuyl and Villarreal in \cite{TV}(also see  in \cite{Z}) showed that the clutter  with the free vertex property is shellable.

\begin{Theorem} \label{Thm1} If a clutter $\mathcal{C}$ with the free vertex property, then $\Delta_{\mathcal{C}}$ is shellable.
\end{Theorem}

\vspace{3mm}The squarefree monomial ideal $I_{m,l,k}$ corresponds to a
clutter (resp. simplicial complex), say
$\mathcal{C}_{m,l,k}$ (resp. $\Delta_{m,l,k}$),
 its edges  are precisely  some such paths of length $m$ in the line graph $L_{n}$. That is,
$E_{\mathcal{C}_{m,l,k}}=\{\{x_{1},\dots, x_{m}\},\{x_{(m-l)+1},
\dots, x_{2(m-l)+l}\},\dots,\\
\{x_{(k-1)(m-l)+1},\dots, x_{k(m-l)+l}\}\}$. Throughout this paper, we will assume that  the clutter with edge set
$E_{\mathcal{C}_{m,l,k}}=\{\{x_{1},\dots, x_{m}\},\{x_{(m-l)+1},\dots, x_{2(m-l)+l}\},
\dots,\{x_{(k-1)(m-l)+1}\\
,\dots,x_{k(m-l)+l}\}\}$
where $l$ is an  integer such that $1\leq l\leq m$. This set  corresponds to a squarefree monomial ideal $I_{m,l,k}$, which is the path ideal of the line graph $L_{n}$, i.e.,
$I_{m,l,k}=I(\mathcal{C}_{m,l,k})$.

\vspace{3mm}
Combining Definition \ref{def2} and Theorem \ref{Thm1}, we then get the following
proposition.

\begin{Proposition} \label{prop1} Let $k,l, m$ be positive integers, $\mathcal{C}_{m,l,k}$ be a clutter with edge set  $E_{\mathcal{C}_{m,l,k}}=\{E_{1},\dots,E_{k}\}$
where $E_{i}=\{x_{(i-1)(m-l)+1},x_{(i-1)(m-l)+2},\dots,
x_{(i-1)(m-l)+m}\}$ for $i=1,\dots,k$ and  $I_{m,l,k}=I(E_{\mathcal{C}_{m,l,k}})$
be the edge ideal of the clutter $\mathcal{C}_{m,l,k}$. Then the quotient ring $R/I_{m,l,k}$ is sequentially Cohen-Macaulay.
\end{Proposition}
 \begin{pf} By theorem \ref{Thm1}, it is enough to prove that $\mathcal{C}_{m,l,k}$
has the free vertex property.
   Let $V=\{x_{i_{1}},\dots,x_{i_{r}},x_{j_{1}},\dots,x_{j_{s}}\}$ be any subset of the set of variables $\{x_{1},\dots,x_{k(m-l)+l}\}$
   and $R'=k[x_{1},\dots,\widehat{x_{i_{1}}},\dots,\widehat{x_{i_{r}}},\dots,\widehat{x_{j_{1}}},
 \dots,\widehat{x_{j_{s}}},
 \dots, x_{k(m-l)+l}]$. \\
 One can assume that
 $I'$ is an ideal of $R'$ minimally generated by monomials $u'_{l_{1}},\dots,u'_{l_{t}}$ with $l_{1}<l_{2}<\cdots<l_{t}$,
 and  $x_{i_{a}}\nmid u'_{l_{b}}$ for any $1\leq a\leq r$, $1\leq b\leq t$, and for any  $1\leq b\leq t$, $u'_{l_{b}}$
 is obtained by dividing
$u_{l_{b}}=\prod\limits_{j=1}^{m}x_{(l_{b}-1)(m-l)+j}$ by the product of all the $x_{j_{c}}$ such that $c\in\{1,\dots,s\}$.
Set $a_{j}=(l_{1}-1)(m-l)+j$ for $j=1,\dots,m$ and   $d=min\{a_{j}\mid j\in \{1,\dots,m\}\  \text{and}\ a_{j}\notin\{j_{1},\dots,j_{s}\}\}$.
It is obvious that $x_{d}$ is a free variable   of $I'$ and   the proof is completed.
\end{pf}

\vspace{3mm}\section{Projective dimension and  regularity of the ideal $I_{m,l,k}$}

\hspace{3mm}In this section, we will provide some  formulas for computing the projective dimension and the regularity of  $I_{m,l,k}$.
As some consequences, we  also give some exact formulas for the depth
of $I_{m,l,k}$.

Let $M$ be a graded $R$-module where $R=K[x_{1},\dots,x_{n}]$.
Associated to $M$ is a minimal graded free resolution of the form

\vspace{3mm}
$0\rightarrow \bigoplus\limits_{j}R(-j)^{\beta_{p,j}(M)}\rightarrow \bigoplus\limits_{j}R(-j)^{\beta_{p-1,j}(M)}\rightarrow \cdots\rightarrow \bigoplus\limits_{j}R(-j)^{\beta_{0,j}(M)}\rightarrow M\rightarrow 0,$
where the maps are exact, $p\leq n$, and $R(-j)$ is the $R$-module obtained by shifting
the degrees of $R$ by $j$. The number
$\beta_{i,j}(M)$, the $(i,j)$-th graded Betti number of $M$,  is
an invariant of $M$ that equals the number of minimal generators of degree $j$ in the
$i$th syzygy module of $M$.
Of particular interest are the following invariants which measure the ¡°size¡± of the minimal graded
free resolution of $I$.
The projective dimension of $I$, denoted pd\,$(I)$, is defined to be
$$\mbox{pd}\,(I):=\mbox{max}\{i\ |\ \beta_{i,j}(I)\neq 0\}.$$
The regularity of $I$, denoted $\mbox{reg}\,(I)$, is defined by
$$\mbox{reg}\,(I):=\mbox{max}\{j-i\ |\ \beta_{i,j}(I)\neq 0\}.$$

\vspace{3mm}We now derive some formulas for $\mbox{pd}\,(I_{m,l,k})$  and $\mbox{reg}\,(I_{m,l,k})$ in some special cases by using some
tools developed in \cite{FHT}.
We let $\mathcal{G}(I)$ denote the unique set of minimal generators of
a monomial ideal $I$.

\begin{Definition} \label{def3}Let $I$  be a monomial ideal, and suppose that there exists monomial
ideals $J$ and $K$ such that $\mathcal{G}(I)$ is the disjoint union of $\mathcal{G}(J)$ and $\mathcal{G}(K)$. Then $I=J+K$
is a Betti splitting if
$$\beta_{i,j}(I)=\beta_{i,j}(J)+\beta_{i,j}(K)+\beta_{i-1,j}(J\cap K)\hspace{2mm}\mbox{for all}\hspace{2mm}i,j\geq 0,$$
where $\beta_{i-1,j}(J\cap K)=0\hspace{2mm}  \mbox{if}\hspace{2mm} i=0$.
\end{Definition}

\vspace{3mm}This formula was first obtained for the total Betti numbers by
Eliahou and Kervaire \cite{EK} and extended to the graded case by Fatabbi \cite{F3}.
In the article \cite{FHT}, the authors describe a number of sufficient conditions for an
ideal $I$ to have a Betti splitting. We shall require the following such condition.

\begin{Theorem}\label{Thm2}(\cite[Corollary 2.7]{FHT}).
Suppose that $I=J+K$ where $\mathcal{G}(J)$ contains all
the generators of $I$ divisible by the variable $x_{i}$ and $\mathcal{G}(K)$ is a nonempty set containing
the remaining generators of $I$. If $J$ has a linear resolution, then $I=J+K$ is a Betti
splitting.
\end{Theorem}

\vspace{3mm}When $I=J+K$ is a Betti
splitting ideal,  Definition \ref{def3} implies the following result:
\begin{Corollary} \label{cor1}
If $I=J+K$ is a Betti splitting, then
\begin{itemize}
 \item[(i)]$\mbox{reg}\,(I)=\mbox{max}\{\mbox{reg}\,(J),\mbox{reg}\,(K),\mbox{reg}\,(J\cap K)-1\}$,
 \item[(ii)] $\mbox{pd}\,(I)=\mbox{max}\{\mbox{pd}\,(J),\mbox{pd}\,(K),\mbox{pd}\,(J\cap K)+1\}$.
\end{itemize}
\end{Corollary}

\vspace{3mm}We need the following Lemma:
\begin{Lemma}
\label{lem1}Let $R_{1}=k[x_{1},\dots,x_{m}]$ and $R_{2}=k[x_{m+1},\dots,x_{n}]$ be two polynomial rings, $I\subseteq R_{1}$ and
$J\subseteq R_{2}$ be two nonzero homogeneous  ideals. Then
\begin{itemize}
\item[(1)]$\mbox{pd}\,(I+J)=\mbox{pd}\,(I)+\mbox{pd}\,(J)+1$,
\item[(2)]$\mbox{reg}\,(I+J)=\mbox{reg}\,(I)+\mbox{reg}\,(J)-1$,
\item[(3)]$\mbox{reg}\,(IJ)=\mbox{reg}\,(I)+\mbox{reg}\,(J)$.
\end{itemize}
\end{Lemma}
\begin{pf}  Let $R=k[x_{1},\dots,x_{n}]$. Then, by Proposition 2.2.20 of \cite{V}, we have that $R/I+J\cong R_{1}/I\otimes_{k} R_{2}/J$.  Hence we get that $\mbox{pd}\,(R/I+J)=\mbox{pd}\,(R_{1}/I)+\mbox{pd}\,(R_{2}/J)$. It follows that
\begin{eqnarray*}\mbox{pd}\,(I+J)&=&\mbox{pd}\,(R/I+J)-1=\mbox{pd}\,(R_{1}/I)+\mbox{pd}\,(R_{2}/J)-1\\
&=&(\mbox{pd}\,(I)+1)+(\mbox{pd}\,(J)+1)-1=\mbox{pd}\,(I)+\mbox{pd}\,(J)+1,
 \end{eqnarray*}

 As for the second and the third assertion,  by Lemma 3.2 of \cite{HT4}, we obtain that $\mbox{reg}\,(R/I+J)=\mbox{reg}\,(R_{1}/I)+\mbox{reg}\,(R_{2}/J)$ and
$\mbox{reg}\,(R/IJ)=\mbox{reg}\,(R_{1}/I)+\mbox{reg}\,(R_{2}/J)+1$.  Therefore, we can conclude that
\begin{eqnarray*}\mbox{reg}\,(I+J)&=&\mbox{reg}\,(R/I+J)+1=\mbox{reg}\,(R_{1}/I)+\mbox{reg}\,(R_{2}/J)+1\\
&=&(\mbox{reg}\,(I)-1)+(\mbox{reg}\,(J)-1)+1=\mbox{reg}\,(I)+\mbox{reg}\,(J)-1,
  \end{eqnarray*}
and
\begin{eqnarray*}\mbox{reg}\,(IJ)&=&\mbox{reg}\,(R/IJ)+1=\mbox{reg}\,(R_{1}/I)+\mbox{reg}\,(R_{2}/J)+2\\
&=&(\mbox{reg}\,(I)-1)+(\mbox{reg}\,(J)-1)+2=\mbox{reg}\,(I)+\mbox{reg}\,(J).
 \end{eqnarray*}
We finished  the proof.
\end{pf}

\vspace{3mm}Now, we prove some main results of this section.

\begin{Theorem}\label{Thm3}
Let $k,l, m, n$ be integers such that $n=k(m-l)+l$ where  $k\geq 1$,  $m\geq 2$ and $l<\lceil \frac{m}{2}\rceil$, here $\lceil \frac{m}{2}\rceil$
denotes the smallest integer $\geq \frac{m}{2}$.
Let $I_{m,l,k}=(u_{1},\dots,u_{k})$ with $u_{i}=\prod\limits_{j=1}^{m}x_{(i-1)(m-l)+j}$ for any $1\leq i\leq k$. Then $\mbox{pd}\,(I_{m,l,k})=k-1$,  $\mbox{reg}\,(I_{m,l,k})=(k-1)(m-l-1)+m$.
\end{Theorem}
\begin{pf} We first claim that  $m-2l-1\geq 0$. In fact, if $m=2s+1$, then $\lceil \frac{m}{2}\rceil=s+1$. By
the  hypothesis, we have that $2l+1\leq 2(\lceil \frac{m}{2}\rceil-1)+1=2s+1=m$. On the other hand,  if $m=2s$,
then $\lceil \frac{m}{2}\rceil=s$. Thus   $2l+1\leq 2(\lceil \frac{m}{2}\rceil-1)+1=2s-1<m$. This proves the claim.
We prove these assertions by induction on $k$. It is clear for $k=1$.
If $k=2$,  we let $J_{1}=I_{m,l,1}$ and $K_{1}=(u_{2})$, which contains all
the generators of $I_{m,l,2}$ divisible by the variable $x_{2m-l}$. Because $K_{1}$ has a linear resolution,
$I_{m,l,2}=J_{1}+K_{1}$ is a Betti
splitting  by Theorem \ref{Thm2} and $J_{1}\cap K_{1}=K_{1}(\prod\limits_{j=1}^{m-l}x_{j})$.  Note that
$J_{1}, K_{1}$ and $J_{1}\cap K_{1}$ are principal  ideals, which implies that $\mbox{pd}\,(J_{1})=\mbox{pd}\,(K_{1})=\mbox{pd}\,(J_{1}\cap K_{1})=0$. Thus, by  Corollary \ref{cor1}, we obtain that $$\mbox{pd}\,(I_{m,l,2})=\mbox{max}\{\mbox{pd}\,(J_{1}),\mbox{pd}\,(K_{1}), \mbox{pd}\,(J_{1}\cap K_{1})+1\}=1.$$
Because  the variables that appear in $K_{1}$ and
 $(\prod\limits_{j=1}^{m-l}x_{j})$ are different,
$\mbox{reg}\,(J_{1}\cap K_{1})=\mbox{reg}\,(J_{1})+\mbox{reg}\,(K_{1})=m+(m-l)$  by Lemma \ref{lem1}. Therefore,
 by Corollary \ref{cor1}, we can conclude that
\begin{eqnarray*}
\mbox{reg}\,(I_{m,l,2})&=&\mbox{max}\{\mbox{reg}\,(J_{1}),\mbox{reg}\,(K_{1}), \mbox{reg}\,(J_{1}\cap K_{1})-1\}\\
&=&\mbox{max}\{m,m,m+(m-l)-1\}=m+(m-l-1).
\end{eqnarray*}
This settles the case $k=2$.

Suppose that $k\geq 3$ and that the statement holds for all $I_{m,l,t}$ with $t<k$. We consider the ideals
$L_{0}=I_{m,l,k}$
and $L_{i}=I_{m,l,k-i-1}+(\prod\limits_{j=1}^{m-l}x_{(k-i-1)(m-l)+j})$ for any $1\leq i\leq k-2$.
We denote  $J_{i}=I_{m,l,k-i}$ for $1\leq i\leq k-1$, $K_{1}=(u_{k})$,
$K_{i}=(\prod\limits_{j=1}^{m-l}x_{(k-i)(m-l)+j})$ for $2\leq i\leq k-1$.
Similar to the case $k=2$, we get that, for $1\leq i\leq k-2$,  $L_{i}=J_{i+1}+K_{i+1}$  is a Betti
splitting. Notice that $J_{i}\cap K_{i}=K_{i}L_{i}$,   for any $1\leq i\leq k-2$,
$J_{k-1}\cap K_{k-1}=K_{k-1}(\prod\limits_{j=1}^{m-l}x_{j})$ and the fact that   the variables that appear in $K_{i}$  and
$L_{i}$ are different and  none of the variables that divide $K_{k-1}$ divide any generator of $\prod\limits_{j=1}^{m-l}x_{j}$,
we obtain that,  for $1\leq i\leq k-2$,
 $\begin{array}{l}
\mbox{pd}\,(J_{i}\cap K_{i})=\mbox{pd}\,(L_{i})=\mbox{max}\{\mbox{pd}\,(J_{i+1}),\mbox{pd}\,(K_{i+1}), \mbox{pd}\,(J_{i+1}\cap K_{i+1})+1\},\\
\mbox{reg}\,(L_{i})=\mbox{max}\{\mbox{reg}\,(J_{i+1}),\mbox{reg}\,(K_{i+1}), \mbox{reg}\,(J_{i+1}\cap K_{i+1})-1\},\\
\mbox{reg}\,(J_{i}\cap K_{i})=\mbox{reg}\,(K_{i}L_{i})=\mbox{reg}\,(K_{i})+\mbox{reg}\,(L_{i})\geq \mbox{reg}\,(K_{i})+1,   \\
\mbox{reg}\,(J_{k-1}\cap K_{k-1})=\mbox{reg}\,(K_{k-1})+\mbox{reg}\,((\prod\limits_{j=1}^{m-l}x_{j}))=2(m-l).
\end{array}~~~(1)$

 Since $J_{k-1}\cap K_{k-1}$ and $K_{i}$ are principal  ideals, $\mbox{pd}\,(J_{k-1}\cap K_{k-1})=\mbox{pd}\,(K_{i})=0$ for $1\leq i\leq k-1$.
By repeated use of  the above equalities (1), the induction assumption $\mbox{pd}\,(J_{i})=k-i-1$, $\mbox{reg}\,(J_{i})=(k-i-1)(m-l-1)+m$  and
  $m-2l-1\geq 0$, we obtain that $\mbox{pd}\,(J_{1}\cap K_{1})=\mbox{pd}\,(L_{1})=k-2$ and $\mbox{reg}\,(J_{1}\cap K_{1})=(k-1)(m-l-1)+m+1$.
It follows that
 \begin{eqnarray*}\mbox{pd}\,(L_{0})\!\!\!&=&\!\!\!\mbox{max}\{\mbox{pd}\,(J_{1}),\mbox{pd}\,(K_{1}), \mbox{pd}\,(J_{1}\cap K_{1})+1\}\\
\!\!\!&=&\!\!\!\mbox{max}\{k-2,0,k-1\}=k-1,\\
\mbox{reg}\,(L_{0})\!\!\!&=&\!\!\!\mbox{max}\{\mbox{reg}\,(J_{1}),\mbox{reg}\,(K_{1}), \mbox{reg}\,(J_{1}\cap K_{1})-1\}\\
\!\!\!&=&\!\!\!\mbox{max}\{(k-2)(m-l-1)+m,m,(k-1)(m-l-1)+m+1-1\}\\
\!\!\!&=&\!\!\!(k-1)(m-l-1)+m,
\end{eqnarray*}
\end{pf}

\vspace{3mm}As a consequence of the above theorem, we have:
\begin{Corollary} \label{cor2}
Let $k,l,m, n$ and $I_{m,l,k}$ be as in Theorem \ref{Thm3},  Then $$\mbox{depth}\,(I_{m,l,k})=n-k+1.$$
\end{Corollary}
\begin{pf} By Auslander-Buchsbaum formula, it follows that $$\mbox{depth}\,(I_{m,l,k})=n-\mbox{pd}\,(I_{m,l,k})=n-k+1.$$
\end{pf}

\vspace{3mm} The following theorem
generalizes Theorem 4.1 of \cite{HT3} and Corollary 4.14 of \cite{AF}.

\begin{Theorem}\label{Thm4}Let $k,l, m, n$ be integers such that $n=k(m-l)+l$ where  $k\geq 1$,  $m\geq 2$ and $\lceil \frac{m}{2}\rceil\leq l<m$.
Let $I_{m,l,k}=(u_{1},\dots,u_{k})$ with $u_{i}=\prod\limits_{j=1}^{m}x_{(i-1)(m-l)+j}$ for any $1\leq i\leq k$.
If $m\equiv 0\,(\mbox{mod}\,(m-l))$
 and we can write $n$ as $n=p(2m-l)+d$ where $0\leq d<2m-l$, then
 \begin{itemize}
\item[(1)]$\mbox{pd}\,(I_{m,l,k})=\left\{\begin{array}{ll}
2p-1& \mbox{if}\ \ d\neq m;\\
2p& \mbox{if}\ \ d=m.
\end{array}\right.$
\item[(2)] $\mbox{reg}\,(I_{m,l,k})=\left\{\begin{array}{ll}
p(2m-l-2)+1& \mbox{if}\ \ d\neq m;\\
p(2m-l-2)+m& \mbox{if}\ \ d=m.
\end{array}\right.$
\end{itemize}
\end{Theorem}
\begin{pf} Let $t=\frac{2m-l}{m-l}$, then $t>2$. In fact, if $t=2$, then $l=0$, contradicting the assumption that $l\geq \lceil \frac{m}{2}\rceil$.
We prove these assertions by induction on $k$.

 The cases $k=1,2$ are from Theorem \ref{Thm3}.
 Suppose that $k\geq 3$ and that the statements hold for all $I_{m,l,s}$ with $s<k$.
 If $3\leq k\leq t$, then $n=(2m-l)+d$ with $d=(k-2)(m-l)<m$. Set $J_{1}=I_{m,l,k-1}$ and $K_{1}=(u_{k})$, we get that
  $J_{1}\cap K_{1}=K_{1}(\prod\limits_{j=1}^{m-l}x_{(k-2)(m-l)+j})$. Thus $\mbox{pd}\,(J_{1}\cap K_{1})=0$ and $\mbox{reg}\,(J_{1}\cap K_{1})=m+(m-l)=2m-l$.
As the number of the variables that appear in $J_{1}$ is $(2m-l)+d-(m-l)$,  using the induction hypothesis, $\mbox{pd}\,(J_{1})=1$ and $\mbox{reg}\,(J_{1})=2m-l-1$. It follows that $\mbox{pd}\,(I_{m,l,k})=\mbox{max}\{\mbox{pd}\,(J_{1}),\mbox{pd}\,(K_{1}), \mbox{pd}\,(J_{1}\cap K_{1})+1\}=1,$
and $\mbox{reg}\,(I_{m,l,k})=\mbox{max}\{\mbox{reg}\,(J_{1}),\mbox{reg}\,(K_{1}),
\mbox{reg}\,(J_{1}\cap K_{1})-1\}=\mbox{max}\{2m-l-1,m,2m-l-1\}=2m-l-1$.
 This proves the assertion for $3\leq k\leq t$.

If $k\geq qt+1$ with $q\geq 1$. Set $J_{1}=I_{m,l,k-1}$ and $K_{1}=(u_{k})$. By similar arguments as in  Theorem \ref{Thm3},
we get that $I_{m,l,k}=J_{1}+K_{1}$  is  a Betti
splitting and  $J_{1}\cap K_{1}=K_{1}(I_{m,l,k-t}+(\prod\limits_{j=1}^{m-l}x_{(k-2)(m-l)+j}))$.
Notice that  the variables that appear in $K_{1}, I_{m,l,k-t}$ and
$(\prod\limits_{j=1}^{m-l}x_{(k-2)(m-l)+j})$ are different,
it follows  that
\begin{eqnarray*}
\mbox{pd}\,(J_{1}\cap K_{1})&=&\mbox{pd}\,(I_{m,l,k-t}+(\prod\limits_{j=1}^{m-l}x_{(k-2)(m-l)+j}))\\
&=&\mbox{pd}\,(I_{m,l,k-t})+\mbox{pd}\,(\prod\limits_{j=1}^{m-l}x_{(k-2)(m-l)+j})+1\\
&=&\mbox{pd}\,(I_{m,l,k-t})+1.
\end{eqnarray*}
where the second equality follows from Lemma \ref{lem1} (1). We distinguish three cases:

(1) If  $k-1=qt$  with $q\geq 1$, then  the numbers of the variables that appear in  $J_{1}$ and $I_{m,l,k-t}$ are $p(2m-l)+l$ and  $(p-1)(2m-l)+m$, respectively. By  inductive assumption, we get that $\mbox{pd}\,(J_{1})=2p-1$, $\mbox{pd}\,(I_{m,l,k-t})=2(p-1)$,
$\mbox{reg}\,(J_{1})=p(2m-l-2)+1$ and $\mbox{reg}\,(I_{m,l,k-t})=(p-1)(2m-l-2)+m$.  Thus
 $$\left.\begin{array}{lll}
\mbox{pd}\,(J_{1}\cap K_{1})&=&\mbox{pd}\,(I_{m,l,k-t})+1=2p-1,\\
\mbox{pd}\,(I_{m,l,k})&=&\mbox{max}\{\mbox{pd}\,(J_{1}),\mbox{pd}\,(K_{1}), \mbox{pd}\,(J_{1}\cap K_{1})+1\}\\
 &=&\mbox{max}\{2p-1,0, (2p-1)+1\}=2p;\\\mbox{reg}\,(J_{1}\cap K_{1})&=&\mbox{reg}\,(K_{1})+\mbox{reg}\,(I_{m,l,k-t}+(\prod\limits_{j=1}^{m-l}x_{(k-2)(m-l)+j}))\\
&=&\mbox{reg}\,(K_{1})+\mbox{reg}\,(I_{m,l,k-t})+\mbox{reg}\,(\prod\limits_{j=1}^{m-l}x_{(k-2)(m-l)+j})-1\\
 &=&m+[(p-1)(2m-l-2)+m]+(m-l)-1\\
  &=&p(2m-l-2)+m+1,\\
\mbox{reg}\,(I_{m,l,k})&=&\mbox{max}\{\mbox{reg}\,(J_{1}),\mbox{reg}\,(K_{1}), \mbox{reg}\,(J_{1}\cap K_{1})-1\}\\
&=&\mbox{max}\{p(2m-l-2)+1, m, p(2m-l-2)+m\}\\
&=&p(2m-l-2)+m.
  \end{array}\right.$$

(2) If $k-1=qt+1$ with $q\geq 1$, then  the numbers of the variables that appear in  $J_{1}$ and $I_{m,l,k-t}$ are $(p-1)(2m-l)+m$ and  $(p-1)(2m-l)$, respectively.
 Thus  by induction, we have that $\mbox{pd}\,(J_{1})=2(p-1)$, $\mbox{pd}\,(I_{m,l,k-t})=2(p-1)-1$,
$\mbox{reg}\,(J_{1})=(p-1)(2m-l-2)+m$ and $\mbox{reg}\,(I_{m,l,k-t})=(p-1)(2m-l-2)+1$.
Therefore, similar to the above assertions, we obtain that
$\mbox{pd}\,(I_{m,l,k})=2p-1$
and $\mbox{reg}\,(I_{m,l,k})=p(2m-l-2)+1.$

(3) If $k-1=qt+c$ with $q\geq 1$ and $2\leq c<t$, then  the numbers of the variables that appear in  $J_{1}$ and $I_{m,l,k-t}$ are
$p(2m-l)+(c-2)(m-l)$ and  $(p-1)(2m-l)+(c-1)(m-l)$, respectively.
 Thus  by induction, we have that $\mbox{pd}\,(J_{1})=2p-1$, $\mbox{pd}\,(I_{m,l,k-t})=2(p-1)-1$,
$\mbox{reg}\,(J_{1})=p(2m-l-2)+1$ and $\mbox{reg}\,(I_{m,l,k-t})=(p-1)(2m-l-2)+1$.
Similarly,  we can conclude that
$\mbox{pd}\,(I_{m,l,k})=2p-1$
and $\mbox{reg}\,(I_{m,l,k})=p(2m-l-2)+1.$
We completed  the proof.
 \end{pf}

\begin{Remark} \label{rem1} Theorem 4.1  of \cite{HT3} and Corollary 4.14 of \cite{AF} are  some  corollaries of  the above theorem  by
specializing to the case that $l=m-1$.
\end{Remark}

\vspace{3mm}As another corollary, we obtain the following result:
\begin{Corollary} \label{cor3}
Let $k,l,m, n$ and $I_{m,l,k}$ be as in Theorem \ref{Thm4},  Then
$$\mbox{depth}\,(I_{m,l,k})=n+2-\lceil \frac{n+(m-l)}{2m-l}\rceil-\lfloor \frac{n+(m-l)}{2m-l}\rfloor.$$
\end{Corollary}
\begin{pf} Let $k-1=qt+c$, where  $q\geq 0$ and $0\leq c<t$. From the proof of the theorem, we get that if $c=0$, then $d=m$, otherwise, $d=(c-1)(m-l)$. Thus, by some straightforward computations, we have that  if $c=0$, then $\lceil \frac{n+(m-l)}{2m-l}\rceil=\lfloor \frac{n+(m-l)}{2m-l}\rfloor=p+1$, otherwise, $\lceil \frac{n+(m-l)}{2m-l}\rceil=p+1$ and
 $\lfloor \frac{n+(m-l)}{2m-l}\rfloor=p$. By Auslander-Buchsbaum formula, we obtain that  $\mbox{depth}\,(I_{m,l,k})=n-\mbox{pd}\,(I_{m,l,k})$, the desired conclusion follows.
\end{pf}

\begin{Theorem}\label{Thm5}Let $k,l, m, n$ be integers such that $n=k(m-l)+l$ where  $k\geq 1$,  $m\geq 2$ and  $\lceil \frac{m}{2}\rceil\leq l<m$.
Let $I_{m,l,k}=(u_{1},\dots,u_{k})$ with $u_{i}=\prod\limits_{j=1}^{m}x_{(i-1)(m-l)+j}$ for any $1\leq i\leq k$.
If $m\equiv s\,(\mbox{mod}\,(m-l))$ with $1\leq s<m-l$
 and we can write $n$ as $n=p(2m-l-s)+d$ where $0\leq d<2m-l-s$, then  $$\mbox{pd}\,(I_{m,l,k})=\left\{\begin{array}{ll}
2p-1& d\neq m;\\
2p& d=m.
\end{array}\right.
$$
\end{Theorem}
\begin{pf} Let $t=\frac{2m-l-s}{m-l}$, then $t>2$ by similar arguments as in Theorem \ref{Thm3}.
We prove these conclusions by induction on $k$.

  The cases $k=1,2$ are from Theorem \ref{Thm3}.
   Suppose that $k\geq 3$ and that the statements hold for all $I_{m,l,s}$ with $s<k$.
 If $3\leq k\leq t$, then $n=(2m-l-s)+d$ with $s+(m-l)\leq d<m$. Set $J_{1}=I_{m,l,k-1}$ and $K_{1}=(u_{k})$, we get that
  $J_{1}\cap K_{1}=K_{1}(\prod\limits_{j=1}^{m-l}x_{(k-2)(m-l)+j})$. Thus $\mbox{pd}\,(J_{1}\cap K_{1})=0$.
As the number of the variables that appear in $J_{1}$ is $(2m-l-s)+d-(m-l)$,  using the induction hypothesis, $\mbox{pd}\,(J_{1})=1$. It follows that $\mbox{pd}\,(I_{m,l,k})=\mbox{max}\{\mbox{pd}\,(J_{1}),\mbox{pd}\,(K_{1}), \mbox{pd}\,(J_{1}\cap K_{1})+1\}=1.$
This proves the assertion for $3\leq k\leq t$.

  If $k\geq qt+1$ with $q\geq 1$. We consider the ideals $L_{0}=I_{m,l,k}$, $J_{1}=I_{m,l,k-1}$, $K_{1}=(u_{k})$, $L_{1}=I_{m,l,k-t}+(\prod\limits_{j=1}^{m-l}x_{(k-2)(m-l)+j})$, $J_{2q}=I_{m,l,k-qt}$,
$K_{2q}=(\prod\limits_{j=1}^{m-l}x_{[k-(q-1)t-2](m-l)+j})$, and for  $1\leq i\leq q-1$,
\begin{center}
 $\begin{array}{l}
 J_{2i}=I_{m,l,k-it}(\Gamma),\\
 J_{2i+1}=I_{m,l,k-it-1}(\Gamma),\\
 K_{2i}=(\prod\limits_{j=1}^{m-l}x_{[k-(i-1)t-2](m-l)+j}),\\
 K_{2i+1}=(\prod\limits_{j=1}^{(t-1)(m-l)}x_{(k-it-1)(m-l)+j}),\\
 L_{2i}=I_{m,l,k-it-1}(\Gamma)+(\prod\limits_{j=1}^{(t-1)(m-l)}x_{(k-it-1)(m-l)+j}),\\ L_{2i+1}=I_{m,l,k-(i+1)t}(\Gamma)+(\prod\limits_{j=1}^{m-l}x_{(k-it-2)(m-l)+j}).
 \end{array}$
 \end{center}
 By similar arguments
as in Theorem \ref{Thm3}, we obtain that, for any $1\leq i\leq 2q$,
 we get that $L_{i}=J_{i+1}+K_{i+1}$  is a Betti splitting
and  $J_{i}\cap K_{i}=K_{i}L_{i}$.  Notice that the variables that appear in $K_{i}$ and
$L_{i}$ are different, we obtain that,  for any  $1\leq i\leq 2q-1$,
$$\mbox{pd}\,(J_{i}\cap K_{i})=\mbox{pd}\,(L_{i})=\mbox{max}\{\mbox{pd}\,(J_{i+1}),\mbox{pd}\,(K_{i+1}), \mbox{pd}\,(J_{i+1}\cap K_{i+1})+1\}.$$
There are three cases to consider:

(1) If $k-1=qt$ for some $q\geq 1$, then $n=k(m-l)+l=(qt+1)(m-l)+l=qt(m-l)+m=q(2m-l-s)+m$. By comparing this with the equality $n=p(2m-l-s)+d$, we have that $q=p$ and $d=m$. The numbers of the variables that appear in $J_{1}$ and $J_{2q}$ are $p(2m-l-s)+l$ and $m$,  respectively.
Similarly, for any $1\leq i\leq q-1$,
the numbers of the variables that appear in $J_{2i+1}$ and $J_{2i}$ are $(p-i)(2m-l-s)+l$ and $(p-i)(2m-l-s)+m$,  respectively.
Hence, by inductive assumption,  $\mbox{pd}\,(J_{1})=2p-1$,  $\mbox{pd}\,(J_{2q})=0$, $\mbox{pd}\,(J_{2i+1})=2(p-i)-1$ and $\mbox{pd}\,(J_{2i})=2(p-i)$ for $1\leq i\leq q-1$.
Note that  $J_{2q}\cap K_{2q}=K_{2q}(\prod\limits_{j=1}^{(t-1)(m-l)}x_{j})$
and $K_{i}$ for $1\leq i\leq 2q$ are principal ideals, we get that  $\mbox{pd}\,(J_{2q}\cap K_{2q})=\mbox{pd}\,(K_{i})=0$. By repeated use of the equality
$\mbox{pd}\,(J_{i}\cap K_{i})=\mbox{max}\{\mbox{pd}\,(J_{i+1}),\mbox{pd}\,(K_{i+1}), \mbox{pd}\,(J_{i+1}\cap K_{i+1})+1\} \ \mbox{for}\ i=2q-1,2q-2,\dots,1,$ we obtain that $\mbox{pd}\,(J_{1}\cap K_{1})=2p-1$. Therefore
\begin{eqnarray*}
\mbox{pd}\,(I_{m,l,k})&=&\mbox{max}\{\mbox{pd}\,(J_{1}),\mbox{pd}\,(K_{1}), \mbox{pd}\,(J_{1}\cap K_{1})+1\}\\
&=&\mbox{max}\{2p-1,0, (2p-1)+1\}=2p.\end{eqnarray*}
This settles the case $k-1=qt$ for some $q\geq 1$.

(2) If $k-1=qt+1$ for some $q\geq 1$, then, similar to the  case (1), we have that
 $q=p+1$ and $d=s$. In this case, the numbers of the variables that appear in $J_{1}$ and $J_{2q}$ are $(p-1)(2m-l-s)+m$ and $1\cdot(2m-l-s)+s$,  respectively. Similarly, for any $1\leq i\leq q-1$, the numbers of the variables that appear in $J_{2i+1}$ and $J_{2i}$ are $(p-i-1)(2m-l-s)+m$ and $(p-i)(2m-l-s)+s$,  respectively.
Hence, by inductive assumption,  $\mbox{pd}\,(J_{1})=2(p-1)$,  $\mbox{pd}\,(J_{2q})=1$, $\mbox{pd}\,(J_{2i+1})=2(p-i-1)$ and $\mbox{pd}\,(J_{2i})=2(p-i)-1$ for $1\leq i\leq q-1$.
Let $L_{2q}=I_{m,l,k-qt-1}+(\prod\limits_{j=1}^{(t-1)(m-l)}x_{(k-qt-1)(m-l)+j})$,
$J_{2q+1}=I_{m,l,k-qt-1}=I_{m,l,1}$,
$K_{2q+1}=(\prod\limits_{j=1}^{(t-1)(m-l)}x_{(m-l)+j})$, then $L_{2q}=J_{2q+1}+K_{2q+1}$  is a Betti splitting and $J_{2q}\cap K_{2q}=K_{2q}L_{2q}$.
Note that   $J_{2q+1}\cap K_{2q+1}=K_{2q+1}(\prod\limits_{j=1}^{m-l}x_{j})$
and $K_{i}$ for $1\leq i\leq 2q+1$ are principal ideals, we get that  $\mbox{pd}\,(J_{2q+1}\cap K_{2q+1})=\mbox{pd}\,(K_{i})=0$. By repeated use of the equality
$\mbox{pd}\,(J_{i}\cap K_{i})=\mbox{max}\{\mbox{pd}\,(J_{i+1}),\mbox{pd}\,(K_{i+1}), \mbox{pd}\,(J_{i+1}\cap K_{i+1})+1\} \ \mbox{for}\ i=2q,2q-1,\dots,1,$ we obtain that $\mbox{pd}\,(J_{1}\cap K_{1})=2(p-1)$. Therefore
\begin{eqnarray*}
\mbox{pd}\,(I_{m,l,k})&=&\mbox{max}\{\mbox{pd}\,(J_{1}),\mbox{pd}\,(K_{1}), \mbox{pd}\,(J_{1}\cap K_{1})+1\}\\
&=&\mbox{max}\{2(p-1),0, 2(p-1)+1\}=2p-1.\end{eqnarray*}
This settles the case $k-1=qt+1$ for some $q\geq 1$.

(3) If $k-1=qt+c$ for some $q\geq 1$ and $2\leq c<t$, then, similar to the  case (1), we have that $p=q+1$ and $d=s+(c-1)(m-l)$.

We claim: $d\neq m$. If $d=m$, then $c-1=\frac{m-s}{m-l}=t-1$. This implies that $c=t$,  contradicting the assumption that
$c<t$. This implies $s+(m-l)\leq d<m-l+(c-1)(m-l)<t(m-l)$.

In this situation, the numbers of the variables that appear in $J_{1}$ and $J_{2q}$ are $p(2m-l-s)+s+(c-2)(m-l)$ and $1\cdot (2m-l-s)+s+(c-1)(m-l)$,  respectively. Similarly, for any $1\leq i\leq q-1$, the numbers of the variables that appear in $J_{2i+1}$ and $J_{2i}$ are $(p-i)(2m-l-s)+s+(c-2)(m-l)$ and $(p-i)(2m-l-s)+s+(c-1)(m-l)$,  respectively.
Hence, by inductive assumption,  $\mbox{pd}\,(J_{1})=2p-1$,  $\mbox{pd}\,(J_{2q})=1$, $\mbox{pd}\,(J_{2i+1})=2(p-i)-1$ and $\mbox{pd}\,(J_{2i})=2(p-i)-1$ for $1\leq i\leq q-1$.
Let $L_{2q}=I_{m,l,k-qt-1}+(\prod\limits_{j=1}^{(t-1)(m-l)}x_{(k-qt-1)(m-l)+j})$,
$J_{2q+1}=I_{m,l,k-qt-1}=I_{m,l,c}$,
$K_{2q+1}=(\prod\limits_{j=1}^{(t-1)(m-l)}x_{c(m-l)+j})$, then  $L_{2q}=J_{2q+1}+K_{2q+1}$  is a Betti splitting, $J_{2q}\cap K_{2q}=K_{2q}L_{2q}$
and  $J_{2q+1}\cap K_{2q+1}=K_{2q+1}(\prod\limits_{j=1}^{m-l}x_{j})$. Similar to the  above case (2), we get that  $\mbox{pd}\,(J_{2q+1}\cap K_{2q+1})=\mbox{pd}\,(K_{i})=0$. By repeated use of the equality
$\mbox{pd}\,(J_{i}\cap K_{i})=\mbox{max}\{\mbox{pd}\,(J_{i+1}),\mbox{pd}\,(K_{i+1}), \mbox{pd}\,(J_{i+1}\cap K_{i+1})+1\} \ \mbox{for}\ i=2q,2q-1,\dots,1,$ we can conclude that $\mbox{pd}\,(J_{1}\cap K_{1})=2(p-1)$. Therefore
\begin{eqnarray*}
\mbox{pd}\,(I_{m,l,k})&=&\mbox{max}\{\mbox{pd}\,(J_{1}),\mbox{pd}\,(K_{1}), \mbox{pd}\,(J_{1}\cap K_{1})+1\}\\
&=&\mbox{max}\{2p-1,0, 2(p-1)+1\}=2p-1.\end{eqnarray*}
The proof is completed.
\end{pf}

\vspace{3mm}An immediate consequence of the above theorem is the following:
\begin{Corollary} \label{cor4}
Let $k,l,m, n,s$ and $I_{m,l,k}$ be as in Theorem \ref{Thm5}.  Then
$$\mbox{depth}\,(I_{m,l,k})=n+2-\lceil \frac{n+m-l-s}{2m-l-s}\rceil-\lfloor \frac{n+m-l-s}{2m-l-s}\rfloor.$$
\end{Corollary}
\begin{pf} Let $k-1=qt+c$, where  $q\geq 0$ and $0\leq c<t$. From the proof of the theorem, we get that if $c=0$, then $d=m$, otherwise, $d=s+(c-1)(m-l)$. Thus, by some straightforward computations, we have that  if $c=0$, then $\lceil \frac{n+(m-l-s)}{2m-l-s}\rceil=\lfloor \frac{n+(m-l-s)}{2m-l-s}\rfloor=p+1$, otherwise, $\lceil \frac{n+(m-l-s)}{2m-l-s}\rceil=p+1$ and
 $\lfloor \frac{n+(m-l-s)}{2m-l-s}\rfloor=p$. By Auslander-Buchsbaum formula, we obtain that  $\mbox{depth}\,(I_{m,l,k})=n-\mbox{pd}\,(I_{m,l,k})$, the desired conclusion follows.
\end{pf}

\vspace{3mm}To conclude, we ask the following open question.
\begin{Problem} Let $k,l,m, n,s$ and $I_{m,l,k}$ be as in Theorem \ref{Thm5}.
Does there exist some methods to compute  the regularity
of  the ideal $I_{m,l,k}$?
\end{Problem}

 \vspace{15mm}

{\bf Acknowledgments}

 \hspace{3mm}
 The author  would like to thank Sara Saeedi Madani who
 point out  some mistakes in my the first version.

\end{document}